\documentclass[12pt, 14paper]{amsart}
\usepackage{amsmath,amsfonts,amssymb,mathtools}

\usepackage{longtable}
\usepackage[mathscr]{eucal}
\usepackage{amsmath, amsthm}
\usepackage{amsbsy}
\usepackage{wasysym}
\usepackage{url}

\theoremstyle{plain}
\usepackage{color}

\newtheorem{theorem}[subsection]{Theorem}

\newtheorem{proposition}[subsection]{Proposition}
\newtheorem{lemma}[subsection]{Lemma}

\newtheorem{example}[subsection]{Example}

\numberwithin{equation}{section}

\input xypic
\xyoption{all}
\usepackage[breaklinks]{hyperref}

\newcommand \ZZ {{\mathbb Z}}

\newcommand \QQ {{\mathbb Q}}

\usepackage{graphicx}
\begin{document}
\title{Class groups of imaginary biquadratic fields}
\author{Kalyan Banerjee,  Kalyan Chakraborty, Arkabrata Ghosh }
\address{SRM University AP}
 \email{kalyan.ba@srmap.edu.in}
 \email{kalyan.c@srmap.edu.in}
  \email{arkabrata.g@srmap.edu.in}
\maketitle

\begin{abstract}
We present two distinct families of imaginary biquadratic fields, each of which contains infinitely many members, with each member having large class groups. Construction of the first family involves elliptic curves and their quadratic twists, whereas to find the other family, we use a combination of elliptic and hyperelliptic curves. Two main results are used, one from Soleng \cite{S} and the other from Banerjee and Hoque \cite{BH}.


\end{abstract}

\section{Introduction}
The class number of number fields is one of the most intriguing and fundamental questions associated with field extensions. 
It is well known that there exist infinitely many quadratic fields, each with a class number divisible by a particular integer. The Cohen-Lenstra heuristic \cite{CL} predicts a certain behavior of the asymptotics of the class number of infinitely many real
quadratic fields.  We refer to \cite{FK06} for additional related work in this direction.

In this paper, we study the class group of imaginary biquadratic fields. Currently, there is only a handful of literature available relating to class groups of biquadratic fields. 
Sime in \cite{Sim} and \cite{Sim1} studied the class group of real biquadratic fields, but very little is known about the class group of imaginary biquadratic fields, and to our knowledge, Athaide et.al.\cite{ACT} is the only work.  

In this direction, we should also mention that the works of Buell \cite{B}, Buell-Call \cite{BC}, Griffin-Ono-Tsai \cite{GOT}, are very interesting and important from the analytic perspective of class groups and its relations with elliptic curves.

This work is motivated towards shedding some light on the class group of imaginary biquadratic fields, and these fields (that we construct) are kind of geometric in the sense that their origin is from elliptic and hyperelliptic curves. The novelty of this work lies that we are extending the Soleng's technique over imaginary quadratic fields and use the result in \cite{BH} to produce element in the class groups of imaginary biquadratic fields.

Let us begin by briefly stating our plan to move from curves to imaginary biquadratic fields. 


We start with an elliptic curve $E$ defined by
\begin{equation}
\label{eq:1.1}
    E: y^2 = x^3 + ax +b,
\end{equation}

where $a,b \in \QQ$. Taking a positive square-free integer $d$ and twisting $E$ over $\QQ(\sqrt{-d})$,
we shall get an elliptic curve,
\begin{equation}
\label{eq:1.2}
E_d: y^2 = x^3 +d^2 ax + d^3 b,
\end{equation}
and it is defined over $\QQ(\sqrt{-d})$ for $E_d$.

 We fix the notation of the twisted elliptic curve as $E_d$ over $\QQ(\sqrt{-d})= K$. In the sequel, $cl(K)$ denotes the class group of $K$.

 Note that if $P = (x,y)$ 
 is a point on $E$, 
then $\bigg((x+s), y\bigg)$ is a point on 
\begin{equation}
\label{eq:1.3}
y^2 = x^3 -3x^2s + (3s^2 +a)x + (-s^3 -as + b)
\end{equation}
for any integer $s$.
A rational point $p=(x,y)=\bigg(A/C^2, B/C^3 \bigg)$ (here $A,B$ and $C$ are integers) gives a proper ideal $ (A, -kB + \sqrt {-s^3 -as + b})$
(where $k \in \ZZ$ and satisfies $kC^3 + s_1 A = 1$ for some $s_1 \in \ZZ$) in 
imaginary quadratic field $\QQ (\sqrt{-s^3 -as + b})$ with discriminant $D= 4( -s^3 -as + b)$ which is not a perfect square.

Thus, taking a torsion point of $E_d$ over $\mathbb{Q}$ and applying Soleng's result (Theorem $4.1$, \cite{S}), we can go from elliptic curves to class groups of imaginary quadratic fields. Then by going up/going down theorem (\cite{DF07}, Chapter 15, \S, 3, Theorem 26), we prove the following result. 
\begin{theorem}
\label{Thm:1.1}
Let $d$ be a square-free integer and let $E_d$ be the quadratic twist of an elliptic curve $E$, defined over $\QQ(\sqrt{-d})$ such that $d^3 b$ is not a square. Let $E_d(\QQ)$ has a primitive element of order $n\leq 16$. Then there are infinitely many biquadratic fields $\QQ(\sqrt{-d},\sqrt{-p^3-d^2pa+d^3b})$ such that the class group of each of them contains an element of order $n$ or $2n$ with $n\neq 2$, provided that $(-p^3-d^2pa+d^3b)$ is not a perfect square.
\end{theorem}

Here primitive torsion point of an elliptic curve is defined in the beginning of \S2.

In the second part, we consider the hyperelliptic curve
\begin{equation}
\label{eq:1.4}
    y^2 =k^2 -x^l
\end{equation}

for large $l ~\text{and}~ k \in \mathbb{Z}$. Then, using the construction described in Theorem $5.1$ (see\cite{BH}), we get an element of order $l$ in the class group of $\QQ(\sqrt{k^2-p^l})$. 
 
 This gives us a family of infinitely many imaginary quadratic fields (there exist infinitely many primes $p$ such that the expression $(k^2-p^l)$ is negative and  is not a square). 
 Thus, by the construction of Soleng (\cite{S}), 
 $\QQ(\sqrt{-q^3-d^2qa+d^3b})$ has an element of order $n$, provided that the elliptic curve $E$ given by \eqref{eq:1.1} has a $n$-torsion (here $n \leq 16$ according to Mazur's theorem). Next, we consider the compositum of the fields 
$$
\QQ(\sqrt{k^2-p^l}) ~\text{and}~ \QQ(\sqrt{-q^3-d^2qa+d^3b})
$$
to obtain a biquadratic field $L$ that has an element of order $ln$. This is possible only if 
$$
 \QQ(\sqrt{k^2-p^l})\cap \QQ(\sqrt{-q^3-d^2qa+d^3b})=\QQ,
$$
that is, they are linearly disjoint. The second result is the following.
\begin{theorem}
\label{Thm:1.2}
 The imaginary biquadratic field $\QQ(\sqrt{k^2-p^l}, \sqrt{-q^3-d^2qa+d^3b})$ has an element of order $ln$ in the class group, where $n\leq 16$ and $n\neq 2$. 
\end{theorem}

In \S 2, we shall briefly discuss Soleng's technique (\cite{S}), which is the stepping stone of this article. Then in \S 3, we shall prove Theorem \eqref{Thm:1.1} followed by an example. Finally, in \S 4, we shall prove Theorem \eqref{Thm:1.2}

\section{Soleng's technique}
Let $P=(x,y)$ be a rational point on the elliptic curve 
\begin{equation}
\label{eq:2.1}
    Y^2=X^3+a_2X^2+a_4X+a_6 
\end{equation}
defined over $\ZZ$. According to Soleng \cite{S}, $P$ is a primitive element if $x,2y$ and $x^2+a_2x+a_4$ are pairwise coprime. By convention, the point at infinity is a primitive point.

The following Proposition (\cite{S}, Proposition 2.1, Page 3) confirms the group structure of these points on the curve.
\begin{proposition}
\label{prop:1}
The set of all primitive points in $E(\QQ)$ is a group. 
\end{proposition}
Let $P=(A/C^2, B/C^3)$ be a point on \eqref{eq:2.1}, then $(A,B)$ is a point on the elliptic curve 
$$
Y^2=X^3+a_2C^2X^2+a_4C^4X+a_6C^6.
$$
Here $a_6\neq 0$ and it is not a perfect square.

Let us consider the integral ideal given by 
$$
I=[A,-B+C^3\sqrt{a_6}].
$$
 Then $I$ defines a proper integral ideal in the order of $\QQ(\sqrt{a_6})$.
The following homomorphism (\cite{S}, Theorem 2.1) gives the desired connection. Let $a_6$ be not a perfect square. Then we have the following theorem.

\begin{theorem}[Theorem 2.1 \cite{S}]\label{thm2}
The homomorphism 
    $$\Phi: E(\QQ)_{prim}\to cl(\QQ[\sqrt{a_6}])$$
is given by 
$$(A/C^2, B/C^3)\to (A,-kB+\sqrt{a_6})$$
where $k$ is an integer satisfying 
$$
kC^3+sA=1
$$
for some $s\in \ZZ$. Here, $cl(\QQ[\sqrt{a_6}])$ denotes the class group of the field $\QQ(\sqrt{a_6})$.
\end{theorem}

Therefore, if we start with an element of order $n$ in $E(\QQ)$ which is a primitive element, then the image of this element under the homomorphism $\Phi$ is of order at most $n$. This homomorphism restricted on non-trivial torsion elements of the group $E(\QQ)_{prim}$  is non-trivial.

\section{The connection between the class groups of biquadratic fields and elliptic curves}
Soleng's main theorem \cite{S}, gives the connection between the rational points of the elliptic curves and the class group of certain number fields.

Now using Theorem (\ref{thm2}), we shall prove Theorem \eqref{Thm:1.1}.

First, we prove the following lemma.

\begin{lemma}
\label{lemma1}
Let $E_d$ be the twist of the elliptic curve $E:y^2=x^3+ax+b$ defined over $\QQ(\sqrt{-d})$,  given by 
$$y^2=x^3+d^2ax+d^3b$$
Then the above homomorphism $\Phi$ is non-trivial from $(E_d)_{prim}$ to 
$$\QQ(\sqrt{-p^3-d^2p+d^3b})$$ for infinitely many primes $p$ (i.e. except for finitely many primes $p$) with loss of primitivity.
\end{lemma}
\begin{proof}
We need to verify that the hypothesis of Theorem \ref{thm2} is satisfied for infinitely many primes $p$ (except for finitely many primes with loss of primitivity). The equation of the twisted elliptic curve is 
$$y^2=x^3+d^2ax+d^3b$$
If $(x,y)$ is a point in this curve then 
$$(x+p,y)$$ is a point on the curve given by the equation
$$y^2=x^3-3px^2+(3p^2+d^2a)x+(-p^3-d^2p+d^3b)$$
Let us call this curve as $E_{1p}$. By a result of\cite{H76}, we know that the number 
$$(-p^3-d^2p+d^3b)$$
is not a perfect square for infinitely many primes $p$ / or except for a finite set of primes. This can also be proven by using a celebrated result of Siegel \cite{Si29}, which states that, any algebraic curve with genus bigger than zero has finitely many integral points on it. Thus  the elliptic curve

\begin{equation}
\label{eq:3.1}
    y^2=x^3+d^2x+d^3b 
\end{equation}

has only finitely many integral solutions. So for any prime $p$ for which the curve \eqref{eq:3.1} has an non integral solution, the expression 
$$(-p^3-d^2p+d^3b)$$
will not a perfect square and these happens except finitely many primes $p$. Hence, in this way we can establish the positive density of the set of primes such that the above number is squarefree.

The other thing we need to check is that if $(x,y)$ is a primitive point on $E_d$, then $(x+p,y)$ is a primitive point on $E_{1p}$. As $(x,y)$ is a primitive point,  we know that $x,2y, x^2+d^2a$ are pairwise coprime. Thus, it follows that $x+p,2y,(x+p)^2+d^2a$ are pairwise coprime for infinitely many primes $p$. Hence $(x+p,y)$ defines a primitive point on $E_{1p}$ for each such prime $p$, and then we can follow Soleng's construction of $\Phi$.
\end{proof}

Now using Theorem (\ref{thm2}) and the above lemma, we shall prove Theorem \eqref{Thm:1.1}.

\begin{proof}
Let us take an elliptic curve $E$ defined over $\QQ$ given by the equation \eqref{eq:1.1}
and we consider its twist over the field $\QQ(\sqrt{-d})$, denoted by $E_d$, which is given by the equation \eqref{eq:1.2}.
This quadratic twist is defined over $\QQ(\sqrt{-d})$.
Now, if we consider the elliptic curve $E_d$ and a non-trivial torsion $\alpha$ on $E_d$, then the specialization of the element $\alpha$ at the point $-p$, denoted by $\alpha_p$ is in the class group of $F=\QQ(\sqrt{-p^3-d^2p+d^3b})$ by the main theorem of Soleng \cite{S} and by \ref{lemma1}. Then by the going up/going down theorem for integral extensions (see \cite{DF07}, Chapter 15, \S 3, Theorem 26), we have an ideal $\beta_p$ in $cl(FK)$ such that $N_{FK/K}(\beta_p)=\alpha_p\neq 0$, and then we have 
$$j_{FK/K}N_{FK/K}(\beta_p)=2\beta_p\;.$$
Here $j_{FK/K},N_{FK/K}$ denote the transfer and the Norm homomorphism at the level of class groups of number fields.

\begin{lemma}
Suppose $\alpha_p$ is of order $n\neq 2$, then we have $\beta_p $ is of order $n$ or $2n$.
\end{lemma}

\begin{proof}
This follows from the identity
$$N_{FK/K}(\beta_p)=\alpha_p$$
and 
$$j_{FK/K}N_{FK/K}(\beta_p)=2\beta_p\;.$$
Therefore, if $n\alpha_p=0$ we have 
$$N_{FK/K}(n\beta_p)=0$$
and hence we have 
$$j_{FK/K}N_{FK/K}(n\beta_p)=2n\beta_p=0$$
So the order of $\beta_p$ divides $2n$. Now suppose that the order is less than $2n$. From the above identity it follows that $n\beta_p$ can be zero if it is not zero then $2n\beta_p=0$. So, in the case when $n\beta_p\neq 0$ we have $2n\beta_p=0$. If the order of $\beta_p$ is less than $2n$ in this case, we have the order divides $2n$ but not $n$. Hence the order is $2$. Then $2\alpha_p=0$ but $n$ is the order of $\alpha_p$ which is not equal to $2$. Therefore we achieve a contradiction. 

In the other case, $n\beta_p=0$, so the order of $\beta_p$ is less than $n$, say $m$. Then by the equality:

$$N_{FK/K}(\beta_p)=\alpha_p$$
we have
$$m\alpha_p=0$$
which is not possible as the order of $\alpha_p$ is $n$ and $m<n$. Therefore the lemma is proven.

\end{proof}

\end{proof}

Now we shall give one, example to illustrate Theorem \eqref{Thm:1.1}


\begin{example}
Let us consider the elliptic curve $E$ given by the equation 
$$E :~ y^2=x^3 + 16$$
over the rational numbers and its twist over 
$$E_d: ~ \QQ(\sqrt{-3})\;.$$
The equation of the twist $E_{-3}$ is given by 
$$y^2=x^3-432\;.$$

Here, $E_{-3}(\mathbb{Q})_{tors} \cong \mathbb{Z}/3$.

 If we specialize the torsion point in $E_{-3}(\QQ)$ generating $\ZZ/3$, on the class group of $L=\QQ(\sqrt{-3}, \sqrt{-p^3-432})$ for a prime $p$,  we get a non-trivial element in the class group of $L$.

 If we call this element $\alpha$, then by the previous theorem \eqref{Thm:1.1}, the order of $\alpha$ is $3$ or $6$.
 
\end{example}

Such examples of elliptic curves can be found in the LMFDB database \cite{LMFDB}, and for each of them, we can make computations related to Soleng's result and produce elements in the class groups of imaginary quadratic fields.


\section{Proof of Theorem ~\texorpdfstring{\eqref{Thm:1.2}}{Thm:1.2}}
\begin{proof}

In the proof of Theorem \eqref{Thm:1.2}, we shall use the techniques of \cite{BH} to prove that $\QQ(\sqrt{k^2-p^l})$ has an element of order $l$ in the class group. Also, by Soleng's technique, $\QQ(\sqrt{-q^3-qd^2a+d^3b})$ has an element of order $n$ where $n\leq 16$ in the class group. Then the compositum has an element of order $ln$ in the class group, provided the intersection of the individual fields in the compositum is $\QQ$.  This is equivalent to $\QQ(\sqrt{k^2-p^l})$ and $ \QQ(\sqrt{-q^3-qd^2a+d^3b})$ are linearly disjoint field extensions $\QQ$. This amounts to showing :

\begin{lemma}The fundamental discriminants satisfy
$$gcd(k^2-p^l, -q^3-qd^2a+d^3b)=1$$ for infinitely many primes $p,q$, such that both of the above numbers are squarefree and congruent to $1$ modulo $4$.
\end{lemma}

\begin{proof}
The claim mentioned above is true  because if for infinitely many primes $p,q$ the above numbers have a common divisor, then we have 
$$d'|k^2-p^l$$
$$d'|-q^3-qd^2a+d^3b$$
and $d'>1$. 
On the other hand, we have $$ -sp^l + t (-q^3) = 1$$
for some integers $s,t$. So, 
$$d'|s(k^2-p^l)+t(-q^3-qd^2a+d^3b)$$
 hence we have 
$$d'|1+sk^2+td^3b-tqd^2a$$
So we have 
$$1+sk^2+td^3b-tqd^2a=md'$$
So, it follows that $d'$ is not a common factor of  $k^2, d^3b-qd^2a$. Since one of the above numbers is $k^2$ which is fixed, then the number of common factors of the two above mentioned numbers is fixed and finite. Hence, the result follows.

\end{proof}



\end{proof}

\section*{Funding and Conflict of Interests/Competing Interests} The author has no financial or non-financial interests to disclose that are directly or indirectly related to the work. The author has no funding sources to report.

\section*{Data availability statement} No outside data was used to prepare this manuscript.

\section*{Acknowledgement}
The authors thank SRM University-AP for providing a suitable atmosphere and support to carry out this research work. The authors thank the anonymous referee for careful reading and suggestions to improve the manuscript.


\end{document}